\documentclass{amsart}
\usepackage[english]{babel}
\usepackage[usenames,dvipsnames]{pstricks}
\usepackage{epsfig}
\usepackage{amssymb,amsmath}
\usepackage[all]{xy}

\newtheorem{theorem}{Theorem}[section]
\newtheorem{proposition}[theorem]{Proposition}
\newtheorem{corollary}[theorem]{Corollary}

\newtheorem{lemma}[theorem]{Lemma}

\newcommand{\TC}{\mathrm{TC}}  
\newcommand{\cat}{\mathrm{cat}} 
\newcommand{\g}{\mathsf{g}}

\newcommand{\ev}{\mathord{\mathrm{ev}}}

\newcommand{\RR}{\mathord{\mathbb{R}}}
\newcommand{\ZZ}{\mathord{\mathbb{Z}}}

\title{Monotonicity of the Schwarz genus}
\author[Petar Pave\v si\'c]{Petar Pave\v si\'c}
\address{University of Ljubljana, Faculty of Mathematics and Physics,  Ljubljana, Slovenia}
\email{\rm{petar.pavesic@fmf.uni-lj.si}}
\begin{document}
\begin{abstract}
\emph{Schwarz genus} $\g(\xi)$ of a fibration $\xi\colon E\to B$ is defined as the minimal integer $n$,  
such that there exists a cover of $B$ by $n$ open sets that admit partial section to $\xi$. Many important concepts,
including Lusternik-Schnirelmann category, Farber's topological complexity and Smale-Vassiliev's complexity of algorithms 
can be naturally expressed as Schwarz genera of suitably chosen fibrations. In this paper we study 
 Schwarz genus 
in relation with certain type of morphisms between fibrations. Our main result is the following: if there exist 
a fibrewise map $f\colon E\to E'$ 
between fibrations $\xi\colon E\to B$ and $\xi'\colon E'\to B$ which induces an $n$-equivalence between 
respective fibres for a sufficiently big $n$, then $\g(\xi)=\g(\xi')$. 
From this we derive several interesting results
relating the topological complexity of a space with the topological complexities of its skeleta and 
subspaces (and similarly for the category). 
For example, we show that if a CW-complexes has high topological complexity (with respect to its dimension and connectivity),
then the topological complexity of its skeleta is an increasing function of the dimension.
\end{abstract}
\maketitle

%============================================================================================================================================================
\section{Introduction}
%============================================================================================================================================================
\label{sec:Introduction}

In this article we will use the standard definitions and results about the Lusternik-Schnirelmann
category and the topological complexity. Interested reader can refer to \cite{CLOT} for 
Lusternik-Schnirelmann category, to \cite{Farber:ITR} for topological complexity of a space, 
and to \cite{Pav:TCMap} for topological complexity of a map. Note however that we follow the 'non-normalized'
convention, so the category and the topological complexity of a contractible space are both equal to 1. 

Our main  objective is to explain the relation between the topological 
complexity of a space and the topological complexity of its skeleta. A well-known result (see \cite[Corollary to Theorem 1]{FHT} or
\cite[Theorem 1.66]{CLOT}) states that if $X$ is a non-contractible CW-complex, then $\cat(X)\ge \cat(X^{(r)})$
where $X^{(r)}$ denotes the $r$-skeleton of $X$. The result is in a sense surprising because 
$\cat(X)$ is a homotopy invariant of $X$, while the homotopy type of the skeleton can vary for different 
CW-decompositions of $X$. As a consequence, the above result restricts the homotopy type of skeleta.
For example, it implies that every skeleton of a (non-contractible) coH-space is itself a coH-space.

The relation between the category of a space and the category of its skeleta does not extended  
directly to topological complexity. For example, the topological complexity is 2 for odd-dimensional spheres
and 3 for even-dimensional spheres (see \cite[Proposition 4.41]{Farber:ITR}). 
Therefore, if we consider the standard CW-decomposition of $S^\infty$, whose skeleta are finite-dimensional
spheres, then the topological complexity of skeleta is an alternating sequence of 2's and 3's, while $\TC(S^\infty)=1$.

In order to understand the causes for different behaviour of two closely related concepts we 
study certain properties of the Schwarz genus of a fibration (see \cite{Schwarz} of \cite[Section 9.3]{CLOT}).
 In fact, the category and the
topological complexity can be both described in terms of the Schwarz genus of 
suitably chosen maps. Section 2 is dedicated to the study of the relations between genera of 
fibrations induced by 
morphisms of fibrations. The main result is Theorem \ref{thm:equality of genera}: it gives sufficient 
conditions on a morphism between fibrations $\xi$ and $\xi'$ over a common base space which imply 
that $\g(\xi)=\g(\xi')$. 
Under similar assumptions but for a morphism between fibrations $\xi$ and $\xi'$ over different 
base spaces we then obtain 
an inequality $\g(\xi)\le\g(\xi')$. Section 3 is split into three subsection in which we apply the general
theory to obtain a series of results that compare topological 
complexity or category of a space to the topological complexity/category of its subspaces.

We will assume throughout this paper that the spaces under consideration are of the homotopy type of 
a CW-complex and have base-points, and that all maps are base-point preserving. Nevertheless we will
systematically omit the base-points from the notation and we will not distinguish notationally between 
a map and its homotopy class.

%============================================================================================================================================================
\section{Comparison of Schwarz genera}
%============================================================================================================================================================
\label{sec:Comparison of Schwarz genera}

Let us recall some basic terminology about fibrations (see \cite{Pav-Pic} for more details).
A (Hurewicz) \emph{fibration} is a triple $(E,\xi,B)$, where the space $B$ is the \emph{base},
$E$ is the \emph{total space} and 
$\xi\colon E\to B$ is a map that has the homotopy lifting property for maps from arbitrary spaces. 
A \emph{morphism} of fibrations is a pair $(f,\bar f)\colon (E,\xi,B)\to (E',\xi',B')$, where 
$f\colon E\to E'$ and $\bar f\colon B\to B'$ are maps, such that the following diagram commutes
$$\xymatrix{
E\ar[r]^f\ar[d]_\xi & E'\ar[d]^{\xi'}\\
B\ar[r]_{\bar f} & B'}
$$
Observe that the map $f$ completely determines the map $\bar f$. 
If $B=B'$ and $\bar f$ is the identity map, we usually contract the notation and write $f$ instead
of $(f,1_B)$. For every $b\in B$ the preimage $\xi^{-1}(b)\subset E$ is called the \emph{fibre} 
of $\xi$ over $b$. In a morphism $(f,\bar f)$ of fibrations the map $f$ clearly sends the fibre 
over $b\in B$ to the fibre over $\bar f(b)\in B'$, so we will occasionally refer to $f$ as a 
\emph{fibrewise} map. A map $\sigma\colon B\to E$ is a \emph{section} of $\xi$ if $\xi\circ\sigma=1_B$.  
More generally, if for some $A\subseteq B$, there exists a map $\sigma\colon A\to E$, such that 
$\xi\circ\sigma$ equals the inclusion of $A$ in $B$, we will call it a \emph{partial section} of $\xi$ over $A$.

In a fibration $(E,\xi,B)$ we will always assume that the base $B$ is path-connected, which in turn implies 
(cf. \cite[Proposition 1.12]{Pav-Pic}) that
all fibres of $\xi$  have the same homotopy type. The fibre over the base-point 
of $B$ will be called \emph{the fibre of $\xi$}. If 
$(f,\bar f)\colon (E,\xi,B)\to (E',\xi',B')$ is a morphism of fibrations, we will denote by
$\tilde f\colon F\to F'$ the induced map between the respective fibres.

Following \cite{Schwarz} we define the \emph{genus} $\g(\xi)$ of the fibration $(E,\xi,B)$
as the minimal integer $n$ for which there exist a cover of $B$ by $n$ 
open sets that admit a partial section of $\xi$. In the context of the Lusternik-Schnirelmann 
category the genus is often called \emph{sectional category} of $\xi$ (see \cite[Section 9.3]{CLOT}).

Observe that genus can be defined for arbitrary maps $\xi\colon E\to B$ by requiring that 
$B$ has a cover by $n$ open sets that admit \emph{homotopy} sections to $\xi$. If $\xi$ is a fibration, then 
every homotopy section can be replaced by a strict section, so the two definitions agree. 
As a matter of fact, most of our results could be easily generalized from fibrations to arbitrary maps.
Our goal in this section is to show that certain kind of morphisms between fibrations induce equality 
between the respective Schwarz genera (see Theorem \ref{thm:equality of genera}). To this end we prove 
several preparatory lemmas.

\begin{lemma}
\label{lem:g(xi)ge g(xi')}
If there exists a morphism  of fibrations $f\colon (E,\xi,B)\to (E',\xi',B)$, then $\g(\xi)\ge \g(\xi')$.
\end{lemma}
\begin{proof}
Consider the following diagram
$$\xymatrix{
  & E \ar[rr]^f \ar[dr]_\xi & & E' \ar[dl]^{\xi'}\\
U \ar@{^(->}[rr] \ar[ur]^\sigma & & B
}$$
If $\sigma$ is a partial section to $\xi$ over $U$, then $f\sigma$ is a partial section to $\xi'$ over $U$, 
and therefore $\g(\xi)\ge \g(\xi')$.
\end{proof}

\begin{lemma}
\label{lem:compatible section}
Let $f\colon (E,\xi,B)\to (E',\xi',B)$ be a morphism of fibrations, and let $\sigma'\colon B\to E'$ 
be a section of $\xi'$. If a map $\sigma\colon B\to E$ satisfies $f\sigma\simeq\sigma'$, then 
$\sigma$ is a homotopy section of $\xi$. 
\end{lemma}
\begin{proof}
Straightforward, since $\xi\sigma=\xi'f\sigma\simeq\xi'\sigma'=1_B$.
\end{proof}

Following \cite[Section VII,6.]{Spanier} we will say that a map $f\colon X\to Y$  is 
an \emph{$n$-equivalence} for $n\ge 1$ if $f$ induces a bijection between the path components 
of $X$ and $Y$, and if the induced homomorphism on the homotopy groups 
$f_\sharp\colon \pi_i(X)\to \pi_i(Y)$ 
is an isomorphism for $0< i<n$ and an epimorphism for $i=n$. A canonical example of an $n$-equivalence is
the inclusion map $X^{(n)}\hookrightarrow X$ of the $n$-skeleton of a CW-complex $X$.
By \cite[Corollary VII,6.23]{Spanier} if $f\colon X\to Y$ is 
an $n$-equivalence, then the induced function $f_*\colon [P,X]\to [P,Y]$ is bijective for every CW-complex
 $P$ of dimension
$\dim(P)\le n-1$, and is surjective if $\dim(P)\le n$.

By analogy, let us say that $f\colon X\to Y$ is a \emph{homology $n$-equivalence} 
if the induced homomorphism on the integral homology groups $f_*\colon H_i(X)\to H_i(Y)$ 
is an isomorphism for $0\le i<n$ and an epimorphism for $i=n$. By \cite[Theorem 7.5.4]{Spanier} an
$n$-equivalence is always a homology $n$-equivalence, and the converse holds if
$X$ and $Y$ are simply-connected.

\begin{lemma}
\label{lem:n-equivalences}
Given a morphism of fibrations $f\colon (E,\xi,B)\to (E',\xi',B)$ the fibrewise map $f\colon E\to E'$ is an $n$-equivalence
if, and only if, the induced map between the respective fibres $\tilde f\colon F\to F'$ is an $n$-equivalence.
\end{lemma}
\begin{proof}
Consider the following 
commutative ladder of exact sequences of homotopy groups of fibrations $(E,\xi,B)$ and $(E',\xi',B)$:
$$\xymatrix{
\cdots \ar[r] & \pi_{i+1}(E) \ar[d]_{f_\sharp} \ar[r]  &  \pi_{i+1}(B) \ar@{=}[d] \ar[r] & 
\pi_i(F) \ar[d]_{\tilde f_\sharp} \ar[r] &  \pi_i(E) \ar[d]_{f_\sharp} \ar[r]  &  
\pi_i(B) \ar@{=}[d] \ar[r] & \cdots                \\
\cdots \ar[r] & \pi_{i+1}(E')                  \ar[r]  &  \pi_{i+1}(B)            \ar[r] &   
\pi_i(F')                         \ar[r] & \pi_i(E')                  \ar[r]  &  
\pi_i(B)            \ar[r] & \cdots
}$$
If $\tilde f_\sharp\colon\pi_i(F)\to \pi_i(F')$ is an isomorphism for $i<n$, then by the Five-lemma
$f_\sharp\colon\pi_i(E)\to \pi_i(E')$ is also an isomorphism for $i<n$. Moreover, if 
$\tilde f_\sharp\colon\pi_{n-1}(F)\to \pi_{n-1}(F')$ is an isomorphism and 
$\tilde f_\sharp\colon\pi_n(F)\to \pi_n(F')$ is an epimorphism, then by the Four-lemma 
$f_\sharp\colon\pi_n(E)\to \pi_n(E')$ is an epimorphism. The converse implication is proved similarly.
\end{proof}

\begin{proposition}
\label{prop:existence of section}
Let $f\colon (E,\xi,B)\to (E',\xi',B)$ be a morphism of fibrations such that the induced map
$\tilde f\colon F\to F'$ is an $n$-equivalence for some $n\ge\dim(B)$.
Then for every section $\sigma'$ of $\xi'$ there exists a compatible section $\sigma$ of $\xi$. 
\end{proposition}
\begin{proof}
By Lemma \ref{lem:compatible section} it is sufficient to find a map $\sigma$ that is a lifting 
of the map $\sigma'$ along $f$ as in the diagram
$$\xymatrix{
  & E \ar[d]^f\\
B \ar@{-->}[ru]^\sigma \ar[r]_{\sigma'} & E'
}$$
By Lemma \ref{lem:n-equivalences} the map $f\colon E\to E'$ is an $n$-equivalence. Since $n\ge\dim(B)$,
the induced function $f_*\colon [B,E]\to [B,E']$ is surjective, which implies that there exists 
$\bar\sigma\colon B\to E$, such that $f\bar\sigma\simeq \sigma'$. By Lemma \ref{lem:compatible section} 
$\bar\sigma$ is 
a homotopy section of $\xi$, and $\bar\sigma$ can be deformed to a strict section $\sigma\colon B\to E$ because
$\xi$ is a fibration. 
\end{proof}

Schwarz \cite[II.1]{Schwarz} introduced a very useful construction that essentially reduces
the computation of the genus to a problem in obstruction theory. Given a fibration $(E.\xi,B)$ 
with fibre $F$ Schwarz defined the \emph{$n$-fold fibrewise join} construction (actually called
\emph{sum} in \cite{Schwarz}) as a fibration $(\ast^n_B E,\xi_n,B)$, where $\ast^n_B E$ is a suitable
subspace of the $n$-fold join $E\ast \cdots\ast E$ , and the
 projection map $\xi_n$ is a fibration whose fibre is the $n$-fold join of fibres $F$, which we 
denote as $\ast^n F$. The main property of Schwarz's construction is stated in the following theorem.

\begin{theorem}[Schwarz, Theorem 3 in \cite{Schwarz}]
\label{thm:Schwarz}
The fibration $(E,\xi,B)$ has genus $\g(\xi)\le n$ if, and only if, the $n$-fold fibrewise 
join $(\ast^n_B E,\xi_n,B)$ admits a section.
\end{theorem}

In other words, $\g(\xi)$ equals the minimal $n$ for which $\xi_n$ admits a section. Observe that 
fibrewise join operation $\ast^n_B$ is functorial, i.e. a morphism 
$$f\colon (E,\xi,B)\to (E',\xi',B)$$ 
induces a morphism 
$$\ast^n_B f\colon (\ast^n_B E,\xi_n,B)\to (\ast^n_B E',\xi'_n,B),$$
whose restriction to the fibres is the usual $n$-fold join of maps 
$\ast^n \tilde f\colon \ast^n F\to \ast^n F'$. 

\begin{lemma}
\label{lem:inductive join}
Let $X$ be a $(c-1)$-connected space, and let the map $f\colon X\to X'$ be a $(c+k)$-equivalence
($c,k\ge 0, c+k\ge 1$). Moreover, let $Y$ be a $(d-1)$-connected space, and let the map $g\colon Y\to Y'$ 
be a $(d+k)$-equivalence ($d\ge 0, d+k\ge 1$). Then $X\ast Y$ is $(c+d)$-connected, and the map
$$f\ast g\colon X\ast Y\to X'\ast Y'$$
is a $(c+d+k+1)$-equivalence.
\end{lemma}
\begin{proof}
Recall the standard homotopy equivalence $X\ast Y\simeq \Sigma (X\wedge Y)$, which together with 
the properties  of the smash product immediately implies that $X\ast Y$ is $(c+d)$-connected. 

Let us first consider the case where at least one of the spaces $X,Y$ is path-connected, i.e. $c+d\ge 1$.
Then 
$X\ast Y$ and $X'\ast Y'$ are 1-connected, so $f\ast g$ is a
$(c+d+k+1)$-equivalence if, and only if, it is a homology 
$(c+d+k+1)$-equivalence. To prove the latter it is sufficient to show that 
$f\wedge g\colon X\wedge Y\to X'\wedge Y'$ is a homology $(c+d+k)$-equivalence.
Consider the morphism of the K\"unneth short exact sequences for the reduced homology of 
the smash product (see \cite[Theorem 5.3.10]{Spanier}):
$$
\xymatrix{
\bigoplus\limits_{i+j=l}(\widetilde H_i(X)\otimes\widetilde H_j(Y))\ \ar@{>->}[r]\ar[d]_{\bigoplus 
(f_*\otimes g_*)} & 
\widetilde H_l(X\wedge Y) \ar@{->>}[r] \ar[d]_{(f\wedge g)_*}& 
\bigoplus\limits_{i+j=l-1}(\widetilde H_i(X)\ast\widetilde H_j(Y)) \ar[d]^{\bigoplus (f_*\ast g_*)}
\\ 
\bigoplus\limits_{i+j=l}(\widetilde H_i(X')\otimes\widetilde H_j(Y'))\ \ar@{>->}[r] & 
\widetilde H_l(X'\wedge Y') \ar@{->>}[r] & \bigoplus\limits_{i+j=l-1}(\widetilde H_i(X')\ast\widetilde H_j(Y'))
}$$
Note that $\widetilde H_i(X)=\widetilde H_i(X')=0$ for $i<c$ and that
$\widetilde H_j(Y)=\widetilde H_j(Y')=0$ for $i<d$. Therefore, if $l< c+d+k$, then  
$\widetilde H_i(X)\otimes\widetilde H_j(Y)\ne 0$ only if  $c\le i< c+k$ and
$d\le j< d+k$. Since $f_*$ and $g_*$ are isomorphisms in that range, we conclude that 
$\bigoplus (f_*\otimes g_*)$ is also an isomorphism. A similar argument shows that 
$\bigoplus (f_*\ast g_*)$ 
is an isomorphism, and therefore the middle map $(f\wedge g)_*$ is also 
an isomorphism. 

If $l=c+d+k$, then we may follow the same line of reasoning to conclude that the first and the last 
summand in 
$\bigoplus (f_*\otimes g_*)$ are epimorphisms, while the remaining summands are isomorphisms, therefore
$\bigoplus (f_*\otimes g_*)$ is an epimorphism. The argument for $\bigoplus (f_*\ast g_*)$ is simpler, 
because $i+j=l-1$, which implies that all summands are isomorphisms and thus  $\bigoplus (f_*\ast g_*)$
is an isomorphism. Four-lemma then implies that $(f\wedge g)_*$ is an epimorphism.

If $c=d=0$ (i.e., $X,X',Y$ and $Y'$ are disconnected), then $k\ge 1$ and $X\ast Y,X'\ast Y'$ are 
path-connected but not necessarily simply-connected. To simplify the notation we assume
that all spaces have at most countably many components, so let $X_0,X_1,\ldots$ be the components
of $X$, and let $Y_0,Y_1,\ldots$ be the components of $Y$. Moreover, assume that the base-point $x_0$
of $X$ is contained in $X_0$, and that the base point $y_0$ of $Y$ is contained in $Y_0$.
By the definition of the smash product we have
$$X\wedge Y=(X_0\wedge Y_0)\  {\scriptstyle\coprod} 
\bigg(\coprod_{i>0} \frac{X_i\times Y_0}{X_i\times y_0}\bigg) {\scriptstyle\coprod}
\bigg(\coprod_{j>0} \frac{X_0\times Y_j}{x_0\times Y_j}\bigg) {\scriptstyle\coprod}
\bigg(\coprod_{i,j>0} X_i\times Y_j\bigg), $$
and there is an analogous description of $X'\wedge Y'$. Since $f\colon X\to X'$ induces a bijection 
between the components of $X$ and $X'$, and similarly for $g\colon Y\to Y'$, it follows that $f\wedge g$ induces 
a bijection between the components of $X\wedge Y$ and $X'\wedge Y'$. By Seifert-van Kampen's theorem, 
$f\ast g\colon X\ast Y\to X'\ast Y'$ induces an isomorphism of respective fundamental groups.
For higher homotopy groups one have to examine the map induced by $f\ast g$ between the respective 
universal covering spaces (or equivalently, consider the induced homomorphisms as homomorphisms
of free $\ZZ G$-modules, where $G=\pi_1(X\ast Y)$). The details are straightforward but tedious so 
we omit them.
\end{proof}

Note that the above lemma essentially states that join operation preserves the property that a map is
equivalence in the $k$ successive dimensions above the connectivity. Thus we have

\begin{proposition}
\label{prop:join}
Assume that $F$ is $(c-1)$-connected ($c\ge 0$) and that $\tilde f\colon F\to F'$ is a $(c+k)$-equivalence for some $k\ge 0$. 
Then $\ast^n \tilde f\colon \ast^n F\to \ast^n F'$ is an $(n(c+1)+k-1)$-equivalence.
\end{proposition}
\begin{proof}
By inductive application of Lemma \ref{lem:inductive join} the $n$-fold join $\ast^n F$ is $(n(c+1)-2)$-connected, hence 
$\pi_i(\ast^n\tilde f)$ is an isomorphism in the following $(k-1)$ dimensions and an epimorphism in dimension $(n(c+1)+k-1)$.
\end{proof}

We are now ready to prove the main theorem of this section.

\begin{theorem}
\label{thm:equality of genera}
Let $(E',\xi',B')$ be a fibration with a $(c-1)$-connected fibre $F'$. If there exists a morphism of fibrations 
$f\colon (E,\xi,B)\to (E',\xi',B')$ such that $\tilde f$ is a $(c+k)$-equivalence for some
$k>\dim(B)-(c+1)\cdot\g(\xi')$, then $\g(\xi)=\g(\xi')$.
\end{theorem}
\begin{proof}
By Lemma \ref{lem:g(xi)ge g(xi')} it is always the case that $\g(\xi)\ge \g(\xi')$, so we 
only need to prove the converse inequality.
By Schwarz's Theorem $n=\g(\xi')$ implies that the fibration $\xi'_n$ admits a global section. 
By naturality of the Schwarz's 
construction there is a morphism of fibrations 
$\ast^n_Bf\colon (\ast^n_B E,\xi_n,B)\to (\ast^n_B E',\xi'_n,B),$
whose restriction to the fibres is an $(n(c+1)+k-1)$-equivalence by Proposition \ref{prop:join}. 
By assumptions of the theorem
$\dim B\le n(c+1)+k-1$, so by Proposition \ref{prop:existence of section} the fibration $\xi_n$ admits 
a section, hence by Schwarz's Theorem $\g(\xi)\le n$.
\end{proof}

We obtain as an immediate consequence the following comparison of genera of fibrations with different 
base spaces.

\begin{corollary}
\label{cor:comparison of genera}
Let $(E',\xi',B')$ be a fibration with a $(c-1)$-connected fibre $F'$. If there exists a morphism 
of fibrations 
$(f,\bar f)\colon (E,\xi,B)\to (E',\xi',B')$ such that $\tilde f$ is a $(c+k)$-equivalence for some
$k>\dim(B)-(c+1)\cdot\g(\xi')$, then $\g(\xi)\le\g(\xi')$.
\end{corollary}
\begin{proof}
We may decompose the morphism $(f,\bar f)$ as in the following diagram, where the middle column represents 
the pullback of $\xi'$ along $\bar f$, and $f$ is equal to the composition $E\to f^*E'\to E'$.
$$\xymatrix{
F  \ar@{^(->}[d] \ar[r]^{\tilde f} & F' \ar@{^(->}[d] \ar@{=}[r] & F'\ar@{^(->}[d]\\
E  \ar[d]_\xi \ar[r] & f^*E' \ar[d]_{f^*\xi'} \ar[r] & E' \ar[d]^{\xi'}\\
B  \ar@{=}[r] & B \ar[r]_{\bar f} & B'
}$$
The fibrewise map $E\to f^*E'$ satisfies the assumptions of Theorem \ref{thm:equality of genera}, therefore 
$\g(\xi)=\g(f^*\xi')$. On the other hand, the genus of a pullback of $\xi'$ is clearly smaller or equal
than the genus of $\xi'$, hence
$\g(\xi)=\g(f^*\xi')\le\g(\xi')$.
\end{proof}

%============================================================================================================================================================
\section{Applications}
%============================================================================================================================================================
\label{sec:Applications}

In this section we are going to compare the values that invariants like the topological complexity 
or category assume on a space $X$  to the values that the same invariants assume on the skeleta
and other subspaces 
of $X$. Recall that the usual approach to the computation of category and topological complexity is to find 
suitable upper and lower estimates. General upper estimates are based on dimension and connectivity
of the spaces involved (see discussion following Theorem \ref{thm:all skeleta above c} for details) while  
the lower estimates usually rely on the multiplicative structures in cohomology. Although the interval between 
the upper and lower estimates can be large, in most cases when the values of topological complexity or 
category are known exactly they are equal (or differ by one) to the general upper bound. Thus, in what follows 
we will normally begin with a general result and then consider the most interesting special case, when 
the value of the invariant is close to the general upper bound. Let us first discuss 
topological complexity of maps.

%--------------------------------------------------------------------------------------------------------------------------------------
\subsection{Topological complexity of maps}\ \\
%--------------------------------------------------------------------------------------------------------------------------------------

Given a continuous map $u\colon X\to Y$ we consider the space $X^I$ of all paths $\alpha\colon I\to X$,
and a map $\xi_u\colon X^I\to X\times Y$ defined as $\xi(\alpha):=(\alpha(0),u(\alpha(1))$. The 
\emph{topological complexity of the map} $u$, denoted $\TC(u)$, is defined as the minimal integer $n$ such that
there exists an increasing sequence of open subsets 
$$\emptyset=U_0\subset U_1\subset\ldots\subset U_n=X\times Y,$$ 
such that each difference $U_i-U_{i-1}$, $i=1,\ldots,n$ admits a continuous partial section 
to the projection $\xi_u$. This concept appeared in  \cite{{Pav:CFKM}} where it was used in order
to measure the manipulation complexity of a robotic device: $X$ and $Y$ were respectively the 
configuration space and the working space of a mechanical system (like a robot arm) 
and $u\colon X\to Y$ was interpreted as the forward kinematic map of the system. The theory was further developed
in \cite{Pav:TCMap}. In particular, if $u$ is a fibration, then by \cite[Lemma 4.1]{Pav:TCMap}
the map $\xi_u\colon X^I\to X\times Y$ is also a fibration. As a consequence, if $u$ is a fibration, 
then its topological complexity can be expressed in terms of the Schwarz's  genus:
$$\TC(u)=\g(\xi_u)$$
(see \cite[Corollary 4.2]{Pav:TCMap}). 

Topological complexity of a single space $X$ is clearly 
equal to the topological complexity of the identity map $1_X$. On the other hand, by 
\cite[Corollary 4.8]{Pav:TCMap} 
the category of $X$ can be retrieved as the topological complexity of the path fibration 
$\ev_1\colon PX\to X$, where $PX$ denotes the space of all based paths in $X$, and $\ev_1$ is the 
evaluation map that sends a based path $\alpha$ to its end-point $\alpha(1)$. 

Theorem 4.9. in \cite{Pav:TCMap} allows further simplification, as it shows that 
$\TC(u)=\g(\eta_u)$, where $\eta_u\colon X\sqcap Y^I\to X\times Y$ is a fibration, whose total space is 
$$X\sqcap Y^I=\{(x,\alpha)\in X\times Y^I\mid f(x)=\alpha(0)\}$$
and $\eta_u(x,\alpha)=(x,\alpha(1))$. It is easy to see that the fibre of $\eta_u$ is the loop space
$\Omega Y$. 

\begin{theorem}
\label{thm:TCmap}
Let $(f,\bar f)\colon (X,u,Y)\to (X',u',Y')$ be a morphism of fibrations, and assume that 
$Y$ is $(c-1)$-connected and that the map $\bar f\colon Y\to Y'$ is an $n$-equivalence, 
for some $n\ge c\ge 1$. 
If $\dim(X\times Y)< c\cdot (\TC(u')-1)+n$, then $\TC(u)\le \TC(u')$.
\end{theorem}
\begin{proof}
First of all, observe that $(f,\bar f)$ determines a morphism of fibration
$$\xymatrix{
X\sqcap Y^I \ar[rr]^{f\times (\bar f\circ -)} \ar[d]_{\eta_u} & & X'\sqcap (Y')^I \ar[d]^{\eta_{u'}}\\
X\times Y \ar[rr]_{f\times \bar f} & & X'\times Y'}
$$
It is easy to check that the induced map between the fibres of $\eta_u$ and $\eta_{u'}$ 
is $\Omega \bar f\colon \Omega Y\to\Omega Y'$. By the assumptions of the theorem 
$\Omega Y$ is $(c-2)$-connected, and $\Omega \bar f$ is an $(n-1)$-equivalence. We may apply 
Corollary \ref{cor:comparison of genera} with $B=X\times Y$ and $k=n-c$ to conclude that 
$c\cdot (\g(\eta_{u'})-1) +n>\dim(X\times Y)$ implies $\g(\eta_{u'})\ge \g(\eta_u)$. The statement of 
the theorem follows immediately from the fact that $\TC(u)=\g(\eta_u)$ and $\TC(u')=\g(\eta_{u'})$.
\end{proof}

%--------------------------------------------------------------------------------------------------------------------------------------
\subsection{Topological complexity of spaces}\ \\
%--------------------------------------------------------------------------------------------------------------------------------------

From the above we can immediately derive a result on topological complexity of spaces as follows. 
Let $X$ be a $(c-1)$-connected
space and let $f\colon X\to X'$ be an $n$-equivalence for some $n\ge c\ge 1$. Then we have a morphism
of trivial fibrations $(f,f)\colon (X,1_X,X)\to (X',1_{X'},X')$ so by Theorem \ref{thm:TCmap}
we obtain
$$\TC(X)\le\TC(X'),$$
provided that $2\dim(X)< c\cdot(\TC(X')-1)+n$. In particular, if $f$ is the inclusion of the $n$-skeleton
$X^{(n)}$ into $X$ we have the following result. 

\begin{corollary}
\label{cor:range}
Let $X$ be a finite-dimensional CW-complex and assume that its $n$-skeleton $X^{(n)}$ is $(c-1)$-connected. 
If $n<c\cdot (\TC(X)-1)$, then $\TC(X)\ge\TC(X^{(n)})$. 
\end{corollary}

If $X\simeq S^m$, then its $n$-skeleton is $(n-1)$-connected 
and at most $n$-dimensional, therefore it is homotopy equivalent to a (possibly empty) wedge of $n$-spheres. 
If $m$ is odd, then $\TC(X)=2$ and the assumptions of the Corollary are not satisfied, 
because the dimension 
of the skeleton cannot be smaller than its connectivity. If $m$ is even, then $\TC(X)=3$ 
and the assumption of the Corollary is that $n<2n$, which holds for all positive $n$. 
Therefore if $X\simeq S^m$ for $m$ even, then $\TC(X^{(n)})\le 3$ for all $n\ge 1$.

By \cite{GLO} we know that $\TC(X)=2$ if, and only if, $X$ is homotopy equivalent to an odd-dimensional sphere. 
Therefore, if $X$ is $(c-1)$-connected and is not homotopy equivalent to a point 
or an odd-dimensional sphere, then
there is always a non-empty range of dimensions, namely $c\le n< c\cdot (\TC(X)-1)$, 
for which $\TC(X)\ge\TC(X^{(n)})$. 
If the topological complexity of $X$ is sufficiently large the above range can cover all skeleta above 
the connectivity of $X$. Thus, by restating the above inequality we obtain the following result.

\begin{theorem}
\label{thm:all skeleta above c}
Let $X$ be a $(c-1)$-connected CW-complex, such that $\TC(X)\ge\frac{\dim(X)}{c}+1$. Then $\TC(X)\ge\TC(X^{(n)})$
for all $n\ge c$.
\end{theorem}

Recall that the topological complexity of $X$ is enclosed between the categories of $X$ and of $X\times X$ 
(see \cite[Proposition 4.19]{Farber:ITR}):
$$\cat(X)\le\TC(X)\le\cat(X\times X),$$
and that for a $(c-1)$-connected CW-complex we have the upper estimate 
$$\cat(X)\le\frac{\dim(X)}{c}+1.$$
(see \cite[Theorem 1.50]{CLOT}).
In particular, if the category of the space equals the above 'dimension-divided-by-connectivity' estimate
(which is often the case), then the assumptions of Theorem \ref{thm:all skeleta above c} are satisfied 
and the topological complexity of $X$ is an upper bound of the topological complexities of its skeleta. 

Furthermore, by combining the above estimates we obtain that 
$$\TC(X)\le\frac{2\dim(X)}{c}+1.$$
Spaces $X$ for which $\TC(X)$ is equal or close to that upper bound were called 
\emph{spaces with high topological complexity} in \cite{FP}. Examples include 
all closed surfaces with the exception of the torus, all complex and quaternionic projective spaces, 
most 3-dimensional lens spaces, configuration spaces, and many other (cf. \cite{Farber:ITR}). Further examples
can be obtained by taking
finite products of the above. For spaces with high topological complexity we can extend 
Theorem \ref{thm:all skeleta above c} to arbitrary subcomplexes.

\begin{corollary}
Let $X$ be a $(c-1)$-connected CW-complex with $\TC(X)\ge\frac{2\dim(X)}{c}$.
Then $\TC(X)\ge\TC(A)$ for every subcomplex $A$ of $X$ containing the $(c+1)$-skeleton of $X$.  
\end{corollary}
\begin{proof}
By cellular approximation theorem the inclusion of $A$ in $X$ is at least a $(c+1)$-equivalence. 
Then $c\cdot(\TC(X)-1)+c+1> 2\dim(A)$ so the statement follows by the discussion preceding 
Corollary \ref{cor:range}.
\end{proof}

Still, we have not been able to rule out possible anomalous behaviour, 
e.g. an increase in topological complexity caused by the removal of a single point 
The following question would be of some interest to applications:

{\bf Question:} Does there exist 
a closed manifold $M$ such that $\TC(M)<\TC(M-x)$ for a single point $x\in M$? 

Recall that the $(k+1)$-skeleton of a space $X$ can be obtained as a mapping cone of a single attaching map
to its $k$-skeleton, therefore $\cat(X^{(k+1)})\le\cat(X^{(k)})+1$. As a consequence the category of 
skeleta cannot 'jump': if $\cat(X)=n$, then every integer $k=1,2,\ldots,n$ must appear 
as the category of some skeleton on $X$. The behaviour 
of the topological complexity can be much more complicated as shown by the sequence of $\TC(\RR P^n)$ 
for $1\le n\le 23$ in \cite[p. 122]{Farber:ITR}. In general we have only a coarse estimate 
\cite[Prop. 4.28]{Farber:ITR}
$$\TC(X^{(k+1)})\le \TC(X^{(k)})+\cat(X^{(k)})+1.$$
Nevertheless, by assuming some control over the category, we are able to show that
the topological complexity is an increasing function along the skeleta.

\begin{proposition}
\label{prop:compare TC}
If $\cat(X)=\dim(X)+1$, then 
$$\TC(X^{(1)})\le \TC(X^{(2)})\le\ldots\le\TC(X).$$
\end{proposition}  
\begin{proof}
Since $\cat(X^{(n)})\le\cat(X^{(n-1)})+1$ for every $n$, it follows that $\cat(X^{(n)})=n+1$ for 
$0\le n\le\dim(X)$. We may thus apply Theorem \ref{thm:all skeleta above c} at each stage of the 
CW-decomposition to show that $\TC(X^{(n)})\ge \TC(X^{(n-1)})$.
\end{proof}

%--------------------------------------------------------------------------------------------------------------------------------------
\subsection{Category of maps and spaces}\ \\
%--------------------------------------------------------------------------------------------------------------------------------------

The \emph{category of a map} $u\colon X\to Y$, denoted $\cat(u)$, is defined as the minimal integer 
$n$ for which there 
exists an open covering $U_1,\ldots,U_n$ of $X$, such that the restrictions 
$f|_{U_i}\colon U_i\to Y$ are null-homotopic (cf. \cite[p. 35]{CLOT}). Clearly, 
the category of a space $X$ is equal to the category of the identity map $1_X$. 
By \cite[Proposition 9.18]{CLOT}
$$\cat(u)=\g(\mu_u)$$
where  $\mu_u\colon X\sqcap PY\to X$ is the pull-back of the path-fibration
$PY\to Y$ along the map $u$, i.e. 
$$X\sqcap PY=\{(x,\alpha)\in X\times PY\mid u(x)=\alpha(1)\}\ \ \text{and}\ \ \mu_u(x,\alpha)=x. $$
Note that we did not require that $u\colon X\to Y$ is a fibration, but it is easy to see that 
the category of a map is equal to the category of its fibrational substitute. As a consequnence, we may 
assume without loss of generality that $u$ is actually a fibration. The proof of the following theorem 
is analogous to that of Theorem \ref{thm:TCmap} so we omit the details.

\begin{theorem}
\label{thm:catmap}
Let $(f,\bar f)\colon (X,u,Y)\to (X',u',Y')$ be a morphism of fibrations, and assume that 
$Y$ is $(c-1)$-connected and that the map $\bar f\colon Y\to Y'$ is an $n$-equivalence, 
for some $n\ge c\ge 1$. 
If $\dim(X)< c\cdot (\cat(u')-1)+n$, then $\cat(u)\le \cat(u')$.
\end{theorem}

If $u$ and $u'$ are taken to be identity maps we obtain a comparison between the category of a space
and the categories of its skeleta and subspaces.
The following theorem is a generalization of  results proved by Felix, Halperin and Thomas \cite{FHT}.

\begin{theorem} Let $X$ be a $(c-1)$-connected CW-complex $X$ ($c\ge 1$).
\begin{itemize}
\item[(1)] 
If $A$ is a subcomplex of $X$ containing $X^{(n)}$ ($n\ge c$)
and such that $$\dim(A)<n+c\cdot(\cat(X)-1), $$ 
then $\cat(A)\le\cat(X)$. 
\item[(2)] 
If $\cat(X)\ge\frac{\dim(X)}{c}$, then $\cat(A)\le\cat(X)$ for every subcomplex $A\le X$
containing the $(c+1)$-skeleton of $X$. 
\item[(3)]
if $X$ is  not contractible,  then $\cat (X^{(n)})\le\cat(X)$ for every $n\ge 1$. 
\item[(4)]
If $X$ is finite-dimensional and $\cat(X)=k$, then 
$$\cat(X^{(\dim(X)-k+1)})\le\ldots\le \cat(X^{(\dim(X)-1)})\le\cat(X).$$
\item[(5)]
If $\cat(X^{(k)})=3$ for some $k$, then 
$$\cat(X^{(k)})\le\cat(X^{(k+1)})\ldots\le \cat(X).$$
\end{itemize}
\end{theorem}
\begin{proof}
If $A$ contains $X^{(n)}$, then the inclusion $A\hookrightarrow X$ is an $n$-equivalence, therefore 
(1) follows directly from Theorem \ref{thm:catmap}. 

By the assumptions of (2), $\cat(X)$ is 
close to the 'dimension-divided-by-connectivity' upper bound for the category. 
The inclusion of $A$ into $X$ 
is a $(c+1)$-equivalence and $c+1+c\cdot(\cat(X)-1)>\dim(X)\ge\dim(A)$, so the statement is an immediate
consequence of (1). Note that the statement remains valid even if $X$ is infinite-dimensional.

If $X$ is not contractible, then $\cat(X)\ge 2$ and so (3) follows from (1). 

For (4) it is sufficient to observe that, like in the proof of Proposition 
\ref{prop:compare TC}, 
if $\cat(X)=k$, then for $n>\dim(X)-k+1$ the $n$-skeleta of $X$ cannot be contractible. Therefore,
we can (3) to $n$-skeleta for $n$ in the stated range.

If $\cat(X^{(n)})\ge 3$, then $\cat(X^{(n+1)})\ge 2$ by \cite[Proposition 2.6]{Cornea}, and then (3) implies
that $\cat(X^{(n+1)})\ge\cat(X^{(n)})\ge 3$. Thus (5) follows by induction.
\end{proof}

\end{document}